\documentclass[10pt]{article}
\usepackage{amsmath}
\usepackage{amssymb}
\usepackage{enumitem}
\usepackage[T1]{fontenc}
\usepackage[utf8]{inputenc}
\usepackage{mathtools}
\usepackage{txfonts}
\usepackage{setspace}
\usepackage{wasysym}
\usepackage[paperheight=27.94cm,paperwidth=21.59cm,left=2.54cm,right=2.54cm,top=2.54cm,bottom=2.54cm]{geometry}
\usepackage[hidelinks]{hyperref}

\begin{document}
	\begin{center}
		{\Large \textbf{NOVEL RESULTS ON SERIES OF FLOOR AND CEILING FUNCTIONS}}
	\end{center}

	\begin{center}
		Dhairya Shah\textsuperscript{1}, Manoj Sahni\textsuperscript{1},Ritu Sahni\textsuperscript{1}, Ernesto Le\'{o}n Castro\textsuperscript{2}, Maricruz Olazabal-Lugo\textsuperscript{3}
	\end{center}

	\begin{center}
		\textit{\textsuperscript{1}Pandit Deendayal Energy University, Gandhinagar, 382007, India.}
	\end{center}

	\begin{center}
		\textit{\textsuperscript{2}Administration Department, Universidad Catolica de la Santisima Concepción,Alonso de Ribera 2850, Concepción, Chile}
	\end{center}
	
	\begin{center}
		\textit{\textsuperscript{3} Universidad Autónoma de Occidente, México}
	\end{center}
	
	\begin{center}
		\href{mailto:dmaths8800@gmail.com}{\textit{dmaths8800@gmail.com}}, \href{mailto:manojsahani117@gmail.com}{\textit{manojsahani117@gmail.com}}, \href{mailto:ritusrivastava1981@gmail.com}{\textit{ritusrivastava1981@gmail.com}}, \href{mailto:eleon@ucsc.cl}{\textit{eleon@ucsc.cl}}
		\href{mailto:maricruz.olazabal@uadeo.mx}{\textit{maricruz.olazabal@uadeo.mx}}
	\end{center}

	ABSTRACT. In the following work, we first propose two (partial summation) formulas involving the floor and ceiling functions. We use principle of mathematical induction to prove the propositions. Another formula relating to the difference of floor and ceiling functions is deduced using aforementioned pair. Finally, in the same section, we propose generalisation of Faulhaber's formula without proof and deduce certain new results using the generalised results. Thereafter, we introduce F$-$Hurwitz and C$-$Hurwitz Zeta functions (infinite series involving floor and ceiling functions respectively) which can be considered as the generalizations of Hurwitz Zeta function. For both infinite series, there exist equivalent series and two distinct methods are used to prove the same. Certain new relations are deduced using new Zeta functions. Thereafter, it is shown that even if new deductions have poles at \( Re\left(s\right) = q\), their differences at the same are convergent. Further some special cases are given for particular values of the Zeta functions. Lastly, certain open problems are provided which might be helpful for further advancements in the field.
	
	\textit{Keywords}: Ceiling function; Floor function; Faulhaber’s formula; Hurwitz Zeta function; Riemann$-$Zeta function
	
		\section{Introduction} \label{sction 1}
	
	Throughout last three centuries, mathematicians such as Gauss, Euler, Bernoulli, Ramanujan, Legendre etc. ventured in the realm of pure mathematics, and especially in the domain of numbers. Some of them discovered very profound formulas, functions and series of their era. The initial development of partial sum of first \( n\) natural numbers \cite{Apostol76} was given by Gauss in the 18th century, and it was generalized for all$-$natural powers of natural numbers, known as Faulhaber's formula (or Bernoulli's formula). In the mid nineteenth century Riemann introduced a working definition for a function of complex variable \cite{MNR17} \( 's' \left(\sum_{n = 1}^{\infty }n^{ - s}\right)\) (A function with interesting history, already introduced for real variables by Euler in the 18th century). The most common generalization of Riemann Zeta function is Hurwitz Zeta function \cite{MNR17}, named after Adolf Hurwitz. Legendre on the other hand introduced the notion of integer part of \( x\) in late 18\textsuperscript{th} century, and Gauss introduced the \(\left[x\right]\) notation for the same, until in the second half of 20\textsuperscript{th} century when the terms floor and ceiling functions were coined.
	
	The authors, reading the history of number theory and studying the behaviour of the Faulhaber’s formula \cite{CG96}, thought on extending the summation to real powers, instead of natural powers. While the work seems interesting, the main problem occurred in exacting the summation and keeping the outcome an integer. As the non$-$integer powers may omit fractional or irrational values, authors had to break the problem into two, having floor and ceiling functions \cite{HW80} of real powers. After achieving the formula, authors considered having the same approach with the Riemann$-$Zeta and Hurwitz$-$Zeta functions \cite{MNR17}. The infinite series with the floor and ceiling function and real powers turned out to be convergent. 
	
	Considering \( p, n\in \mathbb{N}, ~ a\in\left(0,1\right], ~ b\in \mathbb{R}^{ + }\) and \(  s , ~ t \in \mathbb{C},\left(Re\left(s\right)>\frac{1}{a}\right),\left(Re\left(t\right)>0\right)\) The following study introduces novel formulas for \( \sum\limits_{i = 1}^{n}\left\lfloor i^{a}\right\rfloor ,\sum\limits_{i = 1}^{n}\left\lceil i^{a}\right\rceil\), \(  \sum\limits_{n = 0}^{\infty }\left(\left\lfloor bn\right\rfloor^{a} + t\right)^{ - s},\sum\limits_{n = 0}^{\infty }\left(\left\lceil bn\right\rceil^{a} + t\right)^{ - s}\) With proof. Whereas,formulas for\( \sum\limits_{i = 1}^{n}\left\lfloor i^{a}\right\rfloor^{p} ,\sum\limits_{i = 1}^{n}\left\lceil i^{a}\right\rceil^{p}\) are proposed without proof. Specific values are taken to comprehend the behaviour these results.
	
	\section{Preliminaries} \label {section 2}
	
	The following results and definitions are useful for our study:
	
	\subsection{Faulhaber’s (Bernoulli's) Formula \cite{CG96} :} \label{subsection 2.1}
	An expression of the sum of the \( q\) powers of the first \( n\) positive integers can be equated as
	
	\begin{equation*}
		\sum_{i = 1}^{n}i^{q} =\frac{1}{q + 1}\ast \sum_{k = 0}^{q}\binom{q + 1}{k}\ast B_{k}\ast n^{q + 1 - k} = \dfrac{B_{q+1}(n+1)-B_{q+1}(0)}{q+1}
	\end{equation*}
	 here \( B_{j}\) is Bernoulli’s number of second kind, defined as  \\ \begin{equation*}
		B_{j} = \sum_{k = 0}^{j}\frac{1}{k + 1}\ast \sum_{t = 0}^{k}\left( - 1\right)^{t}\ast\left(t + 1\right)^{j}\ast\begin{pmatrix}
			k \\ 
			t \\ 
		\end{pmatrix}
	\end{equation*}
	
	And \(B_{n}(x)\) is Bernoulli Polynomial of order \(n\).
		
	For example, if \( q = 2\) and as \( B_{0} = 1, B_{1} =\frac{1}{2},B_{2} =\frac{1}{6}\)
	
	\begin{equation*}
		\sum_{i = 1}^{n}i^{2} =\frac{1}{3}\ast \sum_{k = 0}^{2}\binom{3}{k}\ast B_{k}\ast n^{3 - k} = \dfrac{B_{3}(n+1)-B_{3}(0)}{3} =\frac{1}{3}n^{3} +\frac{1}{2}n^{2} +\frac{1}{6}n = \dfrac{n(n+1)(2n+1)}{6}
	\end{equation*}

	\subsection{Riemann$-$Zeta function \cite{MNR17}: } The Riemann$-$Zeta function  \( \zeta\left(s\right)\) is a function of a complex variable \( s\) defined as infinite sum
	
	\begin{equation*}
		\zeta\left(s\right) = \sum_{n = 1}^{\infty }\frac{1}{n^{s}} =\frac{1}{\Gamma\left(s\right)}\int_{0}^{\infty }\frac{x^{s - 1}}{e^{x} - 1}dx,  where \Gamma\left(s\right) = \int_{0}^{\infty }t^{s - 1}e^{ - t}dt
	\end{equation*}

	The function converges for all complex value of \( s\) when \( Re\left(s\right)>1\) and defines 
	
	\begin{equation*}
		\zeta\left(s\right) =\frac{1}{1^{s}} +\frac{1}{2^{s}} +\frac{1}{3^{s}} + \cdots 
	\end{equation*}

	\subsection{Hurwitz$-$Zeta function \cite{MNR17}: } The Hurwitz$-$Zeta function  \( \zeta\left(s,t\right)\) is a function of a complex variables \( s\) and \( t\)defined as infinite sum
	
	\begin{equation*}
		\zeta\left(s,t\right) = \sum_{n = 0}^{\infty }\frac{1}{\left(n + t\right)^{s}} =\frac{1}{\Gamma\left(s\right)}\int_{0}^{\infty }\frac{x^{s - 1}e^{ - tx}}{1 - e^{x}}dx,  where ~ \Gamma\left(s\right) = \int_{0}^{\infty }x^{s - 1}e^{ - x}dx
	\end{equation*}

	The series is absolutely convergent for all complex value of \( s\) and \( t\) when \( Re\left(s\right)>1\) and \( Re\left(t\right)>0\).
	
	\subsection{Floor and Ceiling functions \cite{HW80}: } \label {subsection 2.4}
	
	The floor function of any real number \(  x\) (denoted by \( \left\lfloor x\right\rfloor\)) gives the greatest integer not greater than \(  x\). For example, \(\left\lfloor 1.4\right\rfloor  = 1,\left\lfloor 2\right\rfloor  = 2,\left\lfloor  - 3.4\right\rfloor  =  - 4\) and \( \left\lfloor  - 2\right\rfloor  =  - 2\).
	
	The ceiling function (denoted by \( \left\lceil x\right\rceil\)), same way gives the smallest integer not smaller than\(  x\). For example, \(\left\lceil 1.4\right\rceil  = 2,\left\lceil 2\right\rceil  = 2,\left\lceil  - 3.4\right\rceil  =  - 3\) and \( \left\lceil  - 2\right\rceil  =  - 2\).
	
	From above we can see that \(\left\lceil x\right\rceil  =\left\lfloor x\right\rfloor  = x\) if and only if \(  x\in \mathbb{Z}\).
	
	\smallskip	
	
	\section{Formulas} \label {section 3}

	\subsection{Results: }
		
	In this section, we propose the results and provide their proofs. Consider \(  n\in \mathbb{N},\) and \(  a\in\left(0,1\right]\) then two new formulas can be defined as follows: \vspace{3mm}
	\\
	\textbf{FORMULA 1:}
	
	\begin{equation*}
		\sum_{i = 1}^{n}\left\lfloor i^{a}\right\rfloor  =\left(n + 1\right)\left\lfloor n^{a}\right\rfloor  - \sum_{t = 1}^{\left\lfloor n^{a}\right\rfloor }\left\lceil t^{\frac{1}{a}}\right\rceil 
	\end{equation*}

	\textbf{FORMULA 2:}
	
	\begin{equation*}
		\sum_{i = 1}^{n}\left\lceil i^{a}\right\rceil  = n\left\lceil n^{a}\right\rceil  +\left\lfloor\left\lceil n^{a}\right\rceil^{\frac{1}{a}}\right\rfloor  - \sum_{t = 1}^{\left\lceil n^{a}\right\rceil }\left\lfloor t^{\frac{1}{a}}\right\rfloor 
	\end{equation*}

	Note: Formula 1 and 2 holds for all \( a\in \mathbb{R}^{ + }\), but for  \( a>1\), value of \(\left\lfloor n^{a}\right\rfloor\) becomes larger than \( n\) and running the summation from 1 to \(\left\lfloor n^{a}\right\rfloor\) becomes less economical, and hence case when \( a>1\) is not considered.
	
	\subsection{Proofs:}
	
	\textbf{FORMULA 1:} Let \( P\left(n\right)\) be the following statement:
	
	\begin{equation} \label{Equation 1}
		P\left(n\right):\sum_{i = 1}^{n}\left\lfloor i^{a}\right\rfloor  =\left(n + 1\right)\left\lfloor n^{a}\right\rfloor  - \sum_{t = 1}^{\left\lfloor n^{a}\right\rfloor }\left\lceil t^{\frac{1}{a}}\right\rceil
	\end{equation}

	\textbf{Basic step of induction}: \( P\left(1\right)\) is clearly true.
	
	\begin{equation*}
		P\left(1\right):  1 =\left(1 + 1\right)1 - 1 = 2\times 1 - 1 = 1\left(\because a\in\left(0,1\right]\Rightarrow \left\lfloor 1^{a}\right\rfloor  =\left\lceil 1^{\frac{1}{a}}\right\rceil  = 1\right)
	\end{equation*}

	\textbf{Inductive step:} Assume \( P\left(k\right)\) is true for some \(  n = k\). One shall prove \( P\left(k\right)\Rightarrow P\left(k + 1\right)\)
	
	\begin{equation}
		P\left(k\right):\sum_{i = 1}^{k}\left\lfloor i^{a}\right\rfloor  =\left(k + 1\right)\left\lfloor k^{a}\right\rfloor  - \sum_{t = 1}^{\left\lfloor k^{a}\right\rfloor }\left\lceil t^{\frac{1}{a}}\right\rceil
	\end{equation}

	It follows that
	
	\begin{equation}
		P\left(k + 1\right): \sum_{i = 1}^{k + 1}\left\lfloor i^{a}\right\rfloor  =\left(k + 1\right)\left\lfloor k^{a}\right\rfloor  - \sum_{t = 1}^{\left\lfloor k^{a}\right\rfloor }\left\lceil t^{\frac{1}{a}}\right\rceil  +\left\lfloor\left(k + 1\right)^{a}\right\rfloor
	\end{equation}

	Consider \( m\in \mathbb{N}\)such that \(  m =\left\lfloor k^{a}\right\rfloor\).
	
	\begin{equation*}
		m =\left\lfloor k^{a}\right\rfloor \Rightarrow m\leq k^{a}<m + 1\Rightarrow  m^{\frac{1}{a}}\leq k<\left(m + 1\right)^{\frac{1}{a}}\left(\because a\in\left(0,1\right]\right)
	\end{equation*}

	\begin{equation*}
		\Rightarrow \left\lceil m^{\frac{1}{a}}\right\rceil \leq k<\left\lceil\left(m + 1\right)^{\frac{1}{a}}\right\rceil \left(\because k\in \mathbb{N}\right)
	\end{equation*}

	\begin{equation*}
		\Rightarrow \left\lceil m^{\frac{1}{a}}\right\rceil \leq k\leq\left\lceil\left(m + 1\right)^{\frac{1}{a}}\right\rceil  - 1\\ 
	\end{equation*}
	This splits up in two cases
	
	(I): \( k =\left\lceil\left(m + 1\right)^{\frac{1}{a}}\right\rceil  - 1\) $\&$ (II):  \(\left\lceil m^{\frac{1}{a}}\right\rceil \leq k<\left\lceil\left(m + 1\right)^{\frac{1}{a}}\right\rceil  - 1\)
	
	Case (I): \( k =\left\lceil\left(m + 1\right)^{\frac{1}{a}}\right\rceil  - 1\Rightarrow k + 1 =\left\lceil\left(m + 1\right)^{\frac{1}{a}}\right\rceil\)
	
	Here, \( k\in \mathbb{N}\Rightarrow k + 1\in \mathbb{N}\Rightarrow\left\lceil\left(m + 1\right)^{\frac{1}{a}}\right\rceil \in \mathbb{N}\Rightarrow\left\lceil\left(m + 1\right)^{\frac{1}{a}}\right\rceil  =\left(m + 1\right)^{\frac{1}{a}}\Rightarrow k + 1 =\left(m + 1\right)^{\frac{1}{a}}\)
	\begin{equation*}
		\Rightarrow\left(k + 1\right)^{a} = m + 1 
	\end{equation*}
	Again \( m\in \mathbb{N}\Rightarrow m + 1\in \mathbb{N}\Rightarrow\left(k + 1\right)^{a}\in \mathbb{N}\Rightarrow\left(k + 1\right)^{a} =\left\lfloor\left(k + 1\right)^{a}\right\rfloor\) (Subsection \ref{subsection 2.4})
	
	\begin{equation*}
		\Rightarrow\left\lfloor\left(k + 1\right)^{a}\right\rfloor  = m + 1 =\left\lfloor k^{a}\right\rfloor  + 1\left(\because m =\left\lfloor k^{a}\right\rfloor\right)
	\end{equation*}

	\( \Rightarrow\left\lfloor\left(k + 1\right)^{a}\right\rfloor  =\left\lfloor k^{a}\right\rfloor  + 1\), if \( k =\left\lceil\left(m + 1\right)^{\frac{1}{a}}\right\rceil  - 1\)
	
	\( \therefore\) From Equation \(\left(3\right)\)
	
	\begin{equation*}
		P\left(k + 1\right): \sum_{i = 1}^{k + 1}\left\lfloor i^{a}\right\rfloor  =\left(k + 1\right)\left\lfloor k^{a}\right\rfloor  - \sum_{t = 1}^{\left\lfloor k^{a}\right\rfloor }\left\lceil t^{\frac{1}{a}}\right\rceil  +\left\lfloor\left(k + 1\right)^{a}\right\rfloor 
	\end{equation*}

	\begin{equation*}
		\Rightarrow P\left(k + 1\right): \sum_{i = 1}^{k + 1}\left\lfloor i^{a}\right\rfloor  =\left(k + 1\right)\left(\left\lfloor\left(k + 1\right)^{a}\right\rfloor  - 1\right) - \sum_{t = 1}^{\left\lfloor\left(k + 1\right)^{a}\right\rfloor  - 1}\left\lceil t^{\frac{1}{a}}\right\rceil  +\left\lfloor\left(k + 1\right)^{a}\right\rfloor 
	\end{equation*}

	\begin{equation*}
		\Rightarrow P\left(k + 1\right): \sum_{i = 1}^{k + 1}\left\lfloor i^{a}\right\rfloor  =\left(k + 1\right)\left\lfloor\left(k + 1\right)^{a}\right\rfloor  +\left\lfloor\left(k + 1\right)^{a}\right\rfloor  -\left(k + 1\right) - \sum_{t = 1}^{\left\lfloor\left(k + 1\right)^{a}\right\rfloor  - 1}\left\lceil t^{\frac{1}{a}}\right\rceil 
	\end{equation*}

	Now, as \( k + 1 =\left\lceil\left(m + 1\right)^{\frac{1}{a}}\right\rceil\) and \( m + 1 =\left\lfloor\left(k + 1\right)^{a}\right\rfloor\) we get \( k + 1 =\left\lceil\left\lfloor\left(k + 1\right)^{a}\right\rfloor^{\frac{1}{a}}\right\rceil\) 
	
	\begin{equation*}
		\Rightarrow P\left(k + 1\right): \sum_{i = 1}^{k + 1}\left\lfloor i^{a}\right\rfloor  =\left(\left(k + 1\right) + 1\right)\left\lfloor\left(k + 1\right)^{a}\right\rfloor  -\left\lceil\left\lfloor\left(k + 1\right)^{a}\right\rfloor^{\frac{1}{a}}\right\rceil  - \sum_{t = 1}^{\left\lfloor\left(k + 1\right)^{a}\right\rfloor  - 1}\left\lceil t^{\frac{1}{a}}\right\rceil 
	\end{equation*}

	\begin{equation*}
		\Rightarrow P\left(k + 1\right): \sum_{i = 1}^{k + 1}\left\lfloor i^{a}\right\rfloor  =\left(\left(k + 1\right) + 1\right)\left\lfloor\left(k + 1\right)^{a}\right\rfloor  - \sum_{t = 1}^{\left\lfloor\left(k + 1\right)^{a}\right\rfloor }\left\lceil t^{\frac{1}{a}}\right\rceil 
	\end{equation*}
	
	\( \therefore P\left(k\right)\Rightarrow P\left(k + 1\right)\)for case (I). \\ 
	Case (II):
	\begin{align*}
		&\left\lceil m^{\frac{1}{a}}\right\rceil \leq k<\left\lceil\left(m + 1\right)^{\frac{1}{a}}\right\rceil  - 1\Rightarrow\left\lceil m^{\frac{1}{a}}\right\rceil  + 1\leq k + 1<\left\lceil\left(m + 1\right)^{\frac{1}{a}}\right\rceil \\ \Rightarrow & \left\lceil m^{\frac{1}{a}}\right\rceil <k + 1<\left\lceil\left(m + 1\right)^{\frac{1}{a}}\right\rceil \Rightarrow m<\left(k + 1\right)^{a}<m + 1 \\ \Rightarrow & \left\lfloor\left(k + 1\right)^{a}\right\rfloor  = m \Rightarrow \left\lfloor\left(k + 1\right)^{a}\right\rfloor  =\left\lfloor k^{a}\right\rfloor  = m
	\end{align*}

	\( \therefore\) From Equation \(\left(3\right)\)
	
	\begin{equation*}
		P\left(k + 1\right): \sum_{i = 1}^{k + 1}\left\lfloor i^{a}\right\rfloor  =\left(k + 1\right)\left\lfloor k^{a}\right\rfloor  - \sum_{t = 1}^{\left\lfloor k^{a}\right\rfloor }\left\lceil t^{\frac{1}{a}}\right\rceil  +\left\lfloor\left(k + 1\right)^{a}\right\rfloor 
	\end{equation*}
	\begin{equation*}
		\Rightarrow P\left(k + 1\right): \sum_{i = 1}^{k + 1}\left\lfloor i^{a}\right\rfloor  =\left(k + 1\right)\left\lfloor\left(k + 1\right)^{a}\right\rfloor  - \sum_{t = 1}^{\left\lfloor\left(k + 1\right)^{a}\right\rfloor }\left\lceil t^{\frac{1}{a}}\right\rceil  +\left\lfloor\left(k + 1\right)^{a}\right\rfloor                   \left(\because\left\lfloor k^{a}\right\rfloor  =\left\lfloor\left(k + 1\right)^{a}\right\rfloor\right)
	\end{equation*}	
	\begin{equation*}
	 \Rightarrow P\left(k + 1\right): \sum_{i = 1}^{k + 1}\left\lfloor i^{a}\right\rfloor  =\left(\left(k + 1\right) + 1\right)\left\lfloor\left(k + 1\right)^{a}\right\rfloor  - \sum_{t = 1}^{\left\lfloor\left(k + 1\right)^{a}\right\rfloor }\left\lceil t^{\frac{1}{a}}\right\rceil 
	\end{equation*}

	\( \therefore P\left(k\right)\Rightarrow P\left(k + 1\right)\)for case (II).

	\begin{equation*}
		\therefore P\left(k\right)\Rightarrow P\left(k + 1\right)    \forall ~ k ~~ such ~ that ~ \left\lceil m^{\frac{1}{a}}\right\rceil \leq k\leq\left\lceil\left(m + 1\right)^{\frac{1}{a}}\right\rceil  - 1
	\end{equation*}

	Hence \( P\left(k + 1\right)\)is true whenever \( P\left(k\right)\) is true.
	
	Hence by Principle of Mathematical Induction \( P\left(n\right)\) is true \( \forall n\in \mathbb{N}.\) \vspace{7mm}
	
	\textbf{FORMULA 2:} Let \( P\left(n\right)\) be the following statement:
	
	\begin{equation} \label{Equation 4}
		P\left(n\right):\sum_{i = 1}^{n}\left\lceil i^{a}\right\rceil  = n\left\lceil n^{a}\right\rceil  +\left\lfloor\left\lceil n^{a}\right\rceil^{\frac{1}{a}}\right\rfloor  - \sum_{t = 1}^{\left\lceil n^{a}\right\rceil }\left\lfloor t^{\frac{1}{a}}\right\rfloor
	\end{equation}

	\textbf{Basic step of induction}: \( P\left(1\right)\) is clearly true.
	
	\begin{equation*}
		P\left(1\right):  1 = 1\left(1\right) + 1 - 1 = 1\left(\because a\in\left(0,1\right]\Rightarrow \left\lfloor 1^{\frac{1}{a}}\right\rfloor  =\left\lceil 1^{a}\right\rceil  =\left\lfloor\left\lceil 1^{a}\right\rceil^{\frac{1}{a}}\right\rfloor  = 1\right)
	\end{equation*}

	\textbf{Inductive step:} Assume \( P\left(k\right)\) is true for some \(  n = k\). One shall prove \( P\left(k\right)\Rightarrow P\left(k + 1\right)\)
	
	\begin{equation}
		P\left(k\right):\sum_{i = 1}^{k}\left\lceil i^{a}\right\rceil  = k\left\lceil k^{a}\right\rceil  +\left\lfloor\left\lceil k^{a}\right\rceil^{\frac{1}{a}}\right\rfloor  - \sum_{t = 1}^{\left\lceil k^{a}\right\rceil }\left\lfloor t^{\frac{1}{a}}\right\rfloor
	\end{equation}

	It follows that
	
	\begin{equation}
		P\left(k + 1\right): \sum_{i = 1}^{k + 1}\left\lceil i^{a}\right\rceil  = k\left\lceil k^{a}\right\rceil  +\left\lfloor\left\lceil k^{a}\right\rceil^{\frac{1}{a}}\right\rfloor  - \sum_{t = 1}^{\left\lceil k^{a}\right\rceil }\left\lfloor t^{\frac{1}{a}}\right\rfloor  +\left\lceil\left(k + 1\right)^{a}\right\rceil
	\end{equation}

	Consider \( m\in \mathbb{N}\)such that \(  m =\left\lceil k^{a}\right\rceil\).
	
	\begin{equation*}
		m =\left\lceil k^{a}\right\rceil \Rightarrow m - 1<k^{a}\leq m\Rightarrow \left(m - 1\right)^{\frac{1}{a}}<k\leq m^{\frac{1}{a}}\left(\because a\in\left(0,1\right]\right)
	\end{equation*}

	\begin{equation*}
		\Rightarrow\left\lfloor\left(m - 1\right)^{\frac{1}{a}}\right\rfloor <k\leq\left\lfloor m^{\frac{1}{a}}\right\rfloor\left(\because k\in \mathbb{N}\right)\\ \Rightarrow\left\lfloor\left(m - 1\right)^{\frac{1}{a}}\right\rfloor  + 1\leq k\leq\left\lfloor m^{\frac{1}{a}}\right\rfloor \\ 
	\end{equation*}
	This splits up in two cases
	
	(I): \( k =\left\lfloor m^{\frac{1}{a}}\right\rfloor\) $\&$ (II):  \(\left\lfloor \left(m - 1\right)^{\frac{1}{a}}\right\rfloor  + 1\leq k<\left\lfloor m^{\frac{1}{a}}\right\rfloor\)
	
	Case (I): \( k =\left\lfloor m^{\frac{1}{a}}\right\rfloor \Rightarrow k = m^{\frac{1}{a}}\left(\because k\in \mathbb{N}\Rightarrow\left\lfloor m^{\frac{1}{a}}\right\rfloor \in \mathbb{N}\Rightarrow\left\lfloor m^{\frac{1}{a}}\right\rfloor  = m^{\frac{1}{a}}\right)\)
	
	\begin{equation*}
		\Rightarrow k + 1 = m^{\frac{1}{a}} + 1\Rightarrow  k + 1>m^{\frac{1}{a}}\Rightarrow\left(k + 1\right)^{a}>m\Rightarrow\left\lceil\left(k + 1\right)^{a}\right\rceil >\left\lceil m\right\rceil 
	\end{equation*}

	\begin{equation*}
		\Rightarrow\left\lceil\left(k + 1\right)^{a}\right\rceil >m\left(\because m\in \mathbb{N}\right)\Rightarrow\left\lceil\left(k + 1\right)^{a}\right\rceil  = m + 1\Rightarrow\left\lceil\left(k + 1\right)^{a}\right\rceil  =\left\lceil k^{a}\right\rceil  + 1\left(\because m =\left\lceil k^{a}\right\rceil\right)
	\end{equation*}

	\( \Rightarrow\left\lceil\left(k + 1\right)^{a}\right\rceil  =\left\lceil k^{a}\right\rceil  + 1\), if \( k =\left\lfloor m^{\frac{1}{a}}\right\rfloor .\)
	
	\( \therefore\) From Equation \(\left(6\right)\)
	
	\begin{equation*}
		P\left(k + 1\right): \sum_{i = 1}^{k + 1}\left\lceil i^{a}\right\rceil  = k\left\lceil k^{a}\right\rceil  +\left\lfloor\left\lceil k^{a}\right\rceil^{\frac{1}{a}}\right\rfloor  - \sum_{t = 1}^{\left\lceil k^{a}\right\rceil }\left\lfloor t^{\frac{1}{a}}\right\rfloor  +\left\lceil\left(k + 1\right)^{a}\right\rceil 
	\end{equation*}

	\begin{equation*}
		\Rightarrow  P\left(k + 1\right): \sum_{i = 1}^{k + 1}\left\lceil i^{a}\right\rceil  = k\left(\left\lceil\left(k + 1\right)^{a}\right\rceil  - 1\right) +\left\lceil\left(k + 1\right)^{a}\right\rceil  +\left\lfloor\left\lceil k^{a}\right\rceil^{\frac{1}{a}}\right\rfloor  - \sum_{t = 1}^{\left\lceil\left(k + 1\right)^{a}\right\rceil  - 1}\left\lfloor t^{\frac{1}{a}}\right\rfloor 
	\end{equation*}

	\begin{equation*}
		\Rightarrow  P\left(k + 1\right): \sum_{i = 1}^{k + 1}\left\lceil i^{a}\right\rceil  =\left(k + 1\right)\left\lceil\left(k + 1\right)^{a}\right\rceil  - k +\left\lfloor\left\lceil k^{a}\right\rceil^{\frac{1}{a}}\right\rfloor  - \sum_{t = 1}^{\left\lceil\left(k + 1\right)^{a}\right\rceil  - 1}\left\lfloor t^{\frac{1}{a}}\right\rfloor 
	\end{equation*}

	(Adding and subtracting \(\left\lfloor\left\lceil\left(k + 1\right)^{a}\right\rceil^{\frac{1}{a}}\right\rfloor )\) Also. for case (I) we have \(\left\lceil k^{a}\right\rceil  = m\) and \(\left\lfloor m^{\frac{1}{a}}\right\rfloor  = k\Rightarrow\left\lfloor\left\lceil k^{a}\right\rceil^{\frac{1}{a}}\right\rfloor  = k\)
	
	\begin{equation*}
		\Rightarrow  P\left(k + 1\right): \sum_{i = 1}^{k + 1}\left\lceil i^{a}\right\rceil  =\left(k + 1\right)\left\lceil\left(k + 1\right)^{a}\right\rceil  - k + k +\left\lfloor\left\lceil\left(k + 1\right)^{a}\right\rceil^{\frac{1}{a}}\right\rfloor  -\left\lfloor\left\lceil\left(k + 1\right)^{a}\right\rceil^{\frac{1}{a}}\right\rfloor  - \sum_{t = 1}^{\left\lceil\left(k + 1\right)^{a}\right\rceil  - 1}\left\lfloor t^{\frac{1}{a}}\right\rfloor 
	\end{equation*}

	\begin{equation*}
		\Rightarrow  P\left(k + 1\right): \sum_{i = 1}^{k + 1}\left\lceil i^{a}\right\rceil  =\left(k + 1\right)\left\lceil\left(k + 1\right)^{a}\right\rceil  +\left\lfloor\left\lceil\left(k + 1\right)^{a}\right\rceil^{\frac{1}{a}}\right\rfloor  -\left\lfloor\left\lceil\left(k + 1\right)^{a}\right\rceil^{\frac{1}{a}}\right\rfloor  - \sum_{t = 1}^{\left\lceil\left(k + 1\right)^{a}\right\rceil  - 1}\left\lfloor t^{\frac{1}{a}}\right\rfloor 
	\end{equation*}

	\begin{equation*}
		\Rightarrow  P\left(k + 1\right): \sum_{i = 1}^{k + 1}\left\lceil i^{a}\right\rceil  =\left(k + 1\right)\left\lceil\left(k + 1\right)^{a}\right\rceil  +\left\lfloor\left\lceil\left(k + 1\right)^{a}\right\rceil^{\frac{1}{a}}\right\rfloor  - \sum_{t = 1}^{\left\lceil\left(k + 1\right)^{a}\right\rceil }\left\lfloor t^{\frac{1}{a}}\right\rfloor 
	\end{equation*}

	\( \therefore P\left(k\right)\Rightarrow P\left(k + 1\right)\) for case (I).
	
	Case (II):  \(\left\lfloor \left(m - 1\right)^{\frac{1}{a}}\right\rfloor  + 1\leq k<\left\lfloor m^{\frac{1}{a}}\right\rfloor \Rightarrow\left\lfloor\left(m - 1\right)^{\frac{1}{a}}\right\rfloor  + 2\leq k + 1<\left\lfloor m^{\frac{1}{a}}\right\rfloor  + 1\)
	
	\begin{equation*}
		\Rightarrow\left\lfloor\left(m - 1\right)^{\frac{1}{a}}\right\rfloor <k + 1\leq\left\lfloor m^{\frac{1}{a}}\right\rfloor \Rightarrow\left(m - 1\right)<\left(k + 1\right)^{a}\leq m\Rightarrow\left\lceil\left(k + 1\right)^{a}\right\rceil  = m
	\end{equation*}

	\begin{equation*}
		\Rightarrow\left\lceil\left(k + 1\right)^{a}\right\rceil  =\left\lceil k^{a}\right\rceil  = m
	\end{equation*}

	\( \therefore\) From Equation \(\left(6\right)\)
	
	\begin{equation*}
		P\left(k + 1\right): \sum_{i = 1}^{k + 1}\left\lceil i^{a}\right\rceil  = k\left\lceil k^{a}\right\rceil  +\left\lfloor\left\lceil k^{a}\right\rceil^{\frac{1}{a}}\right\rfloor  - \sum_{t = 1}^{\left\lceil k^{a}\right\rceil }\left\lfloor t^{\frac{1}{a}}\right\rfloor  +\left\lceil\left(k + 1\right)^{a}\right\rceil 
	\end{equation*}

	\begin{equation*}
		\Rightarrow  P\left(k + 1\right): \sum_{i = 1}^{k + 1}\left\lceil i^{a}\right\rceil  = k\left\lceil\left(k + 1\right)^{a}\right\rceil  +\left\lfloor\left\lceil\left(k + 1\right)^{a}\right\rceil^{\frac{1}{a}}\right\rfloor  - \sum_{t = 1}^{\left\lceil\left(k + 1\right)^{a}\right\rceil }\left\lfloor t^{\frac{1}{a}}\right\rfloor  +\left\lceil\left(k + 1\right)^{a}\right\rceil 
	\end{equation*}

	\begin{equation*}
		\Rightarrow  P\left(k + 1\right): \sum_{i = 1}^{k + 1}\left\lceil i^{a}\right\rceil  =\left(k + 1\right)\left\lceil\left(k + 1\right)^{a}\right\rceil  +\left\lfloor\left\lceil\left(k + 1\right)^{a}\right\rceil^{\frac{1}{a}}\right\rfloor  - \sum_{t = 1}^{\left\lceil\left(k + 1\right)^{a}\right\rceil }\left\lfloor t^{\frac{1}{a}}\right\rfloor 
	\end{equation*}

	\( \therefore P\left(k\right)\Rightarrow P\left(k + 1\right)\) for case (II).
	
	\begin{equation*}
		\therefore P\left(k\right)\Rightarrow P\left(k + 1\right)        \forall ~ k ~ ~such ~ that~ \left\lfloor\left(m - 1\right)^{\frac{1}{a}}\right\rfloor  + 1\leq k\leq\left\lfloor m^{\frac{1}{a}}\right\rfloor 
	\end{equation*}

	Hence \( P\left(k + 1\right)\) is true whenever \( P\left(k\right)\) is true.
	
	Hence by Principle of Mathematical Induction \( P\left(n\right)\) is true \( \forall n\in \mathbb{N}.\)
	
	\subsection{Deductions}
		
		\textbf{DEDUCTION 1:} Subtracting formula 1 from formula 2 with some basic modification we get:
		
		\begin{equation*}
			\sum_{i = 1}^{n}\left(1 -\left(\left\lceil i^{a}\right\rceil  -\left\lfloor i^{a}\right\rfloor\right)\right) = \sum_{t = 1}^{\left\lfloor n^{a}\right\rfloor }\left(1 -\left(\left\lceil t^{\frac{1}{a}}\right\rceil  -\left\lfloor t^{\frac{1}{a}}\right\rfloor\right)\right), a\in\left(0,1\right]
		\end{equation*}

		We are considering \(  a =\frac{1}{q},q\in \mathbb{N}\) for the following:
		
		Under same condition on \( a\), deduction 1 reduces to: \\ \begin{equation*}
			\sum_{i = 1}^{n}\left(\left\lceil i^{\frac{1}{q}}\right\rceil  -\left\lfloor i^{\frac{1}{q}}\right\rfloor\right) = n -\left\lfloor n^{\frac{1}{q}}\right\rfloor 
		\end{equation*}
		
		\textbf{DEDUCTION 2: }Formula 1 reduces to:
		
		\begin{equation*}
			\sum_{i = 1}^{n}\left\lfloor i^{\frac{1}{q}}\right\rfloor  =\left(n + 1\right)\left\lfloor n^{\frac{1}{q}}\right\rfloor  -\frac{B_{q+1}\left(\left\lfloor n^{\frac{1}{q}}\right\rfloor+1\right)-B_{q+1}(0)}{q+1}
		\end{equation*}

		\textbf{DEDUCTION 3:} Formula 2 reduces to:
		
		\begin{equation*}
			\sum_{i = 1}^{n}\left\lceil i^{\frac{1}{q}}\right\rceil  =\left(n\right)\left\lceil n^{\frac{1}{q}}\right\rceil  +\left\lceil n^{\frac{1}{q}}\right\rceil^{q} -\frac{B_{q+1}\left(\left\lceil n^{\frac{1}{q}}\right\rceil +1\right)-B_{q+1}(0)}{q+1}
		\end{equation*}

		\textbf{Proofs: } Deductions needs no proof, as 1\textsuperscript{st} is simply subtraction of Formula 1 and 2. And 2\textsuperscript{nd} and 3\textsuperscript{rd} are simply the case when \( a\) is of form \(\frac{1}{q}\), proofs of which can be simply be understood using the proofs of the formulas. \vspace{5mm}
		
		{\large \textbf{Results for Specific Values}} \cite{GKP89}
		
		\textbf{(I): }Taking \(  a = 1\) in formula 1 and 2, both of them reduces to the Gauss formula. (Take formula 1 for example)
		
	\begin{equation*}
		\sum_{i = 1}^{n}\left\lfloor i^{\frac{1}{1}}\right\rfloor  =\left(n + 1\right)\left\lfloor n^{\frac{1}{1}}\right\rfloor  - \sum_{i = 1}^{\left\lfloor n^{\frac{1}{1}}\right\rfloor }i^{1}\Rightarrow \sum_{i = 1}^{n}i =\left(n + 1\right)n - \sum_{i = 1}^{n}i\left(\because\left\lfloor x^{1}\right\rfloor  = x ,\forall x\in \mathbb{N}\right)
	\end{equation*}

	\begin{equation*}
		\Rightarrow 2\sum_{i = 1}^{n}i =\left(n + 1\right)n\Rightarrow \sum_{i = 1}^{n}i =\frac{\left(n + 1\right)n}{2}
	\end{equation*}

	For \( a =\frac{1}{2}\) simplifying the RHS:\textbf{ }
	
	\vspace{1\baselineskip}
	\textbf{(I):} Formula 1 becomes:
	
	\begin{equation*}
		\sum_{i = 1}^{n}\left\lfloor\sqrt{i}\right\rfloor  =\frac{6n\left\lfloor\sqrt{n}\right\rfloor - 2\left\lfloor\sqrt{n}\right\rfloor^{3} -  3\left\lfloor\sqrt{n}\right\rfloor^{2} +  5\left\lfloor\sqrt{n}\right\rfloor }{6}
	\end{equation*}

	\textbf{(II):} Formula 2 becomes:
	
	\begin{equation*}
		\sum_{i = 1}^{n}\left\lceil\sqrt{i}\right\rceil  =\frac{6n\left\lceil\sqrt{n}\right\rceil - 2\left\lceil\sqrt{n}\right\rceil^{3} +  3\left\lceil\sqrt{n}\right\rceil^{2} -\left\lceil\sqrt{n}\right\rceil }{6}
	\end{equation*}

	All three of which are available in known literature. With these formulas one can go for \( a =\frac{1}{3},\frac{1}{4},\frac{1}{5},\ldots\) .
	
	\subsection{Generalisation of Faulhaber's Formula:}
	
	In this section we propose generalisations of Faulhaber's formula in two parts using Floor and Ceiling functions. \\
	\textbf{Floor Function Based:}
	
	\[_{}^{F}B_{p}^{a}(n) = \sum_{i = 1}^{n}\left\lfloor i^{a} \right\rfloor^{p} = \sum_{t = 1}^{\left\lfloor n^{a} \right\rfloor}{t^{p}\left( \left\lceil \left( t + 1 \right)^{\frac{1}{a}} \right\rceil - \left\lceil t^{\frac{1}{a}} \right\rceil \right)} + \left\lfloor n^{a} \right\rfloor^{p}\left( n - \left\lceil \left( \left\lfloor n^{a} \right\rfloor + 1 \right)^{\frac{1}{a}} \right\rceil + 1 \right)\]

	\textbf{Deduction A:}

\begin{align*}
	_{}^{F}B_{p}^{\frac{1}{q}}(n) = & \sum_{i = 1}^{n}\left\lfloor \sqrt[q]{i} \right\rfloor^{p} = \sum_{i = 1}^{\lfloor \sqrt[q]{n}\rfloor}{\sum_{t = 0}^{q - 1} \begin{pmatrix}
			q \\ 
			t \\ 
		\end{pmatrix} i^{p + t}} + \left\lfloor \sqrt[q]{n} \right\rfloor^{p}\left( n - \left( \left\lfloor \sqrt[q]{n} \right\rfloor + 1 \right)^{q} + 1 \right) \\ 
	\Rightarrow ~ _{}^{F}B_{p}^{\frac{1}{q}}(n) =  &  \sum_{t = 0}^{q - 1} \begin{pmatrix}
		q \\ 
		t \\ 
	\end{pmatrix}\left\{ \frac{B_{p + t + 1}\left( \left\lfloor \sqrt[q]{n} \right\rfloor + 1 \right) - B_{p + t + 1}\left( 0 \right)}{p + t + 1} \right\} + \left\lfloor \sqrt[q]{n} \right\rfloor^{p}\left( n - \left( \left\lfloor \sqrt[q]{n} \right\rfloor + 1 \right)^{q} + 1 \right)
	\end{align*}
	
	\textbf{Ceiling Function Based:}
	
	\begin{align*}
		_{}^{C}B_{p}^{a}(n) = \sum_{i = 1}^{n}\left\lceil i^{a} \right\rceil^{p} = \sum_{t = 1}^{\left\lceil n^{a} \right\rceil}{t^{p}\left( \left\lfloor t^{\frac{1}{a}} \right\rfloor - \left\lfloor \left( t - 1 \right)^{\frac{1}{a}} \right\rfloor \right)} - \left( \left\lfloor \left\lceil n^{a} \right\rceil^{\frac{1}{a}} \right\rfloor - n \right)\left\lceil n^{a} \right\rceil^{p}\
	\end{align*}

	\textbf{Deduction B:}
	
	\begin{align*}
		& _{}^{C}B_{p}^{\frac{1}{q}}(n) = \sum_{i = 1}^{n}\left\lceil \sqrt[q]{i} \right\rceil^{p} = \sum_{i = 1}^{\left\lceil \sqrt[q]{n} \right\rceil}{\sum_{t = 0}^{q - 1}{\left( - 1 \right)^{q - t + 1}\begin{pmatrix}
		q \\ 
		t \\ 
		\end{pmatrix}}i^{p + t}} - \left( \left\lceil \sqrt[q]{n} \right\rceil^{q} - n \right)\left\lceil \sqrt[q]{n} \right\rceil^{p}\\ \Rightarrow ~ & _{}^{C}B_{p}^{\frac{1}{q}}(n) = \sum_{t = 0}^{q - 1}{\left( - 1 \right)^{q - t + 1}\begin{pmatrix}
		q \\ 
		t \\ 
		\end{pmatrix}}\left\{ \frac{B_{p + t + 1}\left( \left\lceil \sqrt[q]{n} \right\rceil + 1 \right) - B_{p + t + 1}\left( 0 \right)}{p + t + 1} \right\} - \left( \left\lceil \sqrt[q]{n} \right\rceil^{q} - n \right)\left\lceil \sqrt[q]{n} \right\rceil^{p}
	\end{align*}
	
	For \(a = 1\) both the partial sums reduce to the Faulhaber's formula, whereas for \(p = 1\) they reduce to equations (\ref{Equation 1}) and (\ref{Equation 4}) respectively.
	
	\section{Series/Zeta Functions}	\label {section 4}			
	In this section, we propose the results and provide their proofs. Consider: \(  n\in \mathbb{N}, a\in\left(0,1\right], b\in \mathbb{R}^{ + }, s\in \mathbb{C}\) and \( t\in \mathbb{C}\)	
	
	\subsection{Results:} \label{subsection 4.1}

	\textbf{SERIES 1 / F$-$Hurwitz Zeta Function:}
	
	F – Hurwitz (Floor) Zeta function is the infinite series given by:
	
	\begin{equation*}
		_{ }^{F}\zeta_{b}^{a}\left(s,t\right) = \sum_{n = 0}^{\infty }\frac{1}{\left(\left\lfloor\left(bn\right)^{a}\right\rfloor  + t\right)^{s}} , Re\left(s\right)>\frac{1}{a} , Re\left(t\right)>0
	\end{equation*}

	However, for the infinite series there exists another infinite series which is equivalent to it:
	
	\begin{equation*}
		\sum_{n = 0}^{\infty }\frac{1}{\left(\left\lfloor\left(bn\right)^{a}\right\rfloor  + t\right)^{s}} =\frac{1}{t^{s}} + \sum_{n = 1}^{\infty }\frac{\left\lceil\frac{\left(n + 1\right)^{\frac{1}{a}}}{b}\right\rceil  -\left\lceil\frac{n^{\frac{1}{a}}}{b}\right\rceil }{\left(n + t\right)^{s}}, Re\left(s\right)>\frac{1}{a} , Re\left(t\right)>0
	\end{equation*}

	\textbf{SERIES 2 / C$-$Hurwitz Zeta Function: }
	
	C – Hurwitz (Ceiling) Zeta function is an infinite series given by:  
	
	\begin{equation*}
		_{  }^{C}\zeta_{b}^{a}\left(s,t\right) = \sum_{n = 0}^{\infty }\frac{1}{\left(\left\lceil\left(bn\right)^{a}\right\rceil  + t\right)^{s}} , Re\left(s\right)>\frac{1}{a} ,Re\left(t\right)>0
	\end{equation*}

	However, for the infinite series there exists another infinite series which is equivalent to it:
	
	\begin{equation*}
		\sum_{n = 0}^{\infty }\frac{1}{\left(\left\lceil\left(bn\right)^{a}\right\rceil  + t\right)^{s}} =\frac{1}{t^{s}} + \sum_{n = 1}^{\infty }\frac{\left\lfloor\frac{n^{\frac{1}{a}}}{b}\right\rfloor  -\left\lfloor\frac{(n - 1)^{\frac{1}{a}}}{b}\right\rfloor }{\left(n + t\right)^{s}}, Re\left(s\right)>\frac{1}{a} ,Re\left(t\right)>0
	\end{equation*}

	\subsection{Proofs and Supplementary Derivations: }
	
	\subsubsection{Proofs}
	
	\textbf{SERIES 1:}
	
	\begin{equation*}
		\sum_{n = 0}^{\infty }\frac{1}{\left(\left\lfloor\left(bn\right)^{a}\right\rfloor  + t\right)^{s}} =\frac{1}{t^{s}} + \sum_{n = 1}^{\infty }\frac{\left\lceil\frac{\left(n + 1\right)^{\frac{1}{a}}}{b}\right\rceil  -\left\lceil\frac{n^{\frac{1}{a}}}{b}\right\rceil }{\left(n + t\right)^{s}} , Re\left(s\right)>\frac{1}{a} , Re\left(t\right)>0
	\end{equation*}

	Let \( m =\left\lfloor\left(bn\right)^{a}\right\rfloor\) for some \( m\in \mathbb{N}\) and let \( f\left(m\right)\) denote function that gives number of consecutive integers,\(  n ,\) for which \( m\) is the particular natural number (i.e., no. of repetition of \(  m\)). Then
	
	\begin{equation*} 
		\sum_{n=0}^{\infty}\frac{1}{\left(\left\lfloor\left(bn\right)^{a}\right\rfloor + t\right)^{s}} =\frac{1}{t^{s}} + \sum_{m = 1}^{\infty }\frac{f\left(m\right)}{\left(m + t\right)^{s}}
	\end{equation*}

	Now,
	
	\begin{equation*}
		m =\left\lfloor\left(bn\right)^{a}\right\rfloor \Rightarrow m\leq\left(bn\right)^{a}<m + 1\Rightarrow  m^{\frac{1}{a}}\leq\left(bn\right)<\left(m + 1\right)^{\frac{1}{a}}\left(\because a\in\left(0,1\right]\right)
	\end{equation*}

	\begin{equation*}
		\Rightarrow \frac{m^{\frac{1}{a}}}{b}\leq n<\frac{\left(m + 1\right)^{\frac{1}{a}}}{b} \Rightarrow \left\lceil\frac{m^{\frac{1}{a}}}{b}\right\rceil \leq n<\left\lceil\frac{\left(m + 1\right)^{\frac{1}{a}}}{b}\right\rceil\left(\because n\in \mathbb{N}\right)
	\end{equation*}

	It follows that \(  n\) is at least \(\left\lceil\frac{m^{\frac{1}{a}}}{b}\right\rceil\) and at most strictly less than \( \left\lceil\frac{\left(m + 1\right)^{\frac{1}{a}}}{b}\right\rceil\)
	
	\( \therefore \) No of consecutive integer \( n\) such that  \( m =\left\lfloor\left(bn\right)^{a}\right\rfloor\) is \(\left\lceil\frac{\left(m + 1\right)^{\frac{1}{a}}}{b}\right\rceil  -\left\lceil\frac{m^{\frac{1}{a}}}{b}\right\rceil\)
	
	\begin{equation*}
		\therefore f\left(m\right) =\left\lceil\frac{\left(m + 1\right)^{\frac{1}{a}}}{b}\right\rceil  -\left\lceil\frac{m^{\frac{1}{a}}}{b}\right\rceil 
	\end{equation*}

	\begin{equation*}
		\therefore \sum_{n = 0}^{\infty }\frac{1}{\left(\left\lfloor\left(bn\right)^{a}\right\rfloor  + t\right)^{s}} =\frac{1}{t^{s}} + \sum_{m = 1}^{\infty }\left\{\frac{\left\lceil\frac{\left(m + 1\right)^{\frac{1}{a}}}{b}\right\rceil  -\left\lceil\frac{m^{\frac{1}{a}}}{b}\right\rceil }{\left(m + t\right)^{s}}\right\}  = \sum_{m = 0}^{\infty }\left\{\frac{\left\lceil\frac{\left(m + 1\right)^{\frac{1}{a}}}{b}\right\rceil  -\left\lceil\frac{m^{\frac{1}{a}}}{b}\right\rceil }{\left(m + t\right)^{s}}\right\} 
	\end{equation*}

	\textbf{SERIES 2:}
	
	\begin{equation*}
		\sum_{n = 0}^{\infty }\frac{1}{\left(\left\lceil\left(bn\right)^{a}\right\rceil  + t\right)^{s}} =\frac{1}{t^{s}} + \sum_{n = 1}^{\infty }\frac{\left\lfloor\frac{n^{\frac{1}{a}}}{b}\right\rfloor  -\left\lfloor\frac{(n - 1)^{\frac{1}{a}}}{b}\right\rfloor }{\left(n + t\right)^{s}}, Re\left(s\right)>\frac{1}{a} ,Re\left(t\right)>0
	\end{equation*}

	Let \( k =\left\lceil\left(bn\right)^{a}\right\rceil\) for some \( k\in \mathbb{N}\) and let \( g\left(k\right)\) denote function that gives number of consecutive integers, \(  n ,\) for which \( k\) is the particular natural number (i.e. no. of repetition of \(  k\)). Then
	
	\begin{equation*}
		\sum_{n = 0}^{\infty }\frac{1}{\left(\left\lceil\left(bn\right)^{a}\right\rceil  + t\right)^{s}} =\frac{1}{t^{s}} + \sum_{k = 1}^{\infty }\frac{g\left(k\right)}{\left(k + t\right)^{s}}
	\end{equation*}

	Now,
	
	\begin{equation*}
		k =\left\lceil\left(bn\right)^{a}\right\rceil \Rightarrow k - 1<\left(bn\right)^{a}\leq k\Rightarrow \left(k - 1\right)^{\frac{1}{a}}<\left(bn\right)\leq k^{\frac{1}{a}}\left(\because a\in\left(0,1\right]\right)
	\end{equation*}

	\begin{equation*}
		\Rightarrow \frac{\left(k - 1\right)^{\frac{1}{a}}}{b}<n\leq\frac{k^{\frac{1}{a}}}{b}\Rightarrow \left\lfloor\frac{(k - 1)^{\frac{1}{a}}}{b}\right\rfloor <n\leq\left\lfloor\frac{k^{\frac{1}{a}}}{b}\right\rfloor\left(\because n\in \mathbb{N}\right)
	\end{equation*}

	It follows that \(  n\) is at least strictly greater than \(\left\lfloor\frac{(k - 1)^{\frac{1}{a}}}{b}\right\rfloor\) and at most \(\left\lfloor\frac{k^{\frac{1}{a}}}{b}\right\rfloor .\)
	
	\( \therefore \)No of consecutive integer \( n\) such that  \( k =\left\lceil\left(bn\right)^{a}\right\rceil\) is \(\left\lfloor\frac{k^{\frac{1}{a}}}{b}\right\rfloor  -\left\lfloor\frac{(k - 1)^{\frac{1}{a}}}{b}\right\rfloor\)
	
	\begin{equation*}
		\therefore g\left(k\right) =\left\lfloor\frac{k^{\frac{1}{a}}}{b}\right\rfloor  -\left\lfloor\frac{(k - 1)^{\frac{1}{a}}}{b}\right\rfloor 
	\end{equation*}

	\begin{equation*}
		\therefore \sum_{n = 0}^{\infty }\frac{1}{\left(\left\lceil\left(bn\right)^{a}\right\rceil  + t\right)^{s}} =\frac{1}{t^{s}} + \sum_{k = 1}^{\infty }\left\{\frac{\left\lfloor\frac{k^{\frac{1}{a}}}{b}\right\rfloor  -\left\lfloor\frac{(k - 1)^{\frac{1}{a}}}{b}\right\rfloor }{\left(k + t\right)^{s}}\right\} 
	\end{equation*}

	\subsubsection{Derivations:}
		
	\textbf{SERIES I:}
	
	Consider the following improper integral (Known as the Gamma Function)
	
	\begin{equation*}
		\Gamma\left(s\right) = \int_{0}^{\infty }y^{s - 1}e^{ - y}dy
	\end{equation*}

	Take \( y = x\left(\left\lfloor\left(bn\right)^{a}\right\rfloor  + t\right)\)
	
	\( \Rightarrow dy = \left(\left\lfloor\left(bn\right)^{a}\right\rfloor  + t\right)dx\), \( x\rightarrow 0\) as \( y\rightarrow 0\) and \( x\rightarrow \infty\) as \( y\rightarrow \infty\)
	
	\begin{equation*}
		\therefore  \Gamma\left(s\right) = \int_{0}^{\infty }\left(x\left(\left\lfloor\left(bn\right)^{a}\right\rfloor  + t\right)\right)^{s - 1}e^{ - x\left(\left\lfloor\left(bn\right)^{a}\right\rfloor  + t\right)}\left(\left\lfloor\left(bn\right)^{a}\right\rfloor  + t\right)dx
	\end{equation*}

	\begin{equation*}
		\therefore  \Gamma\left(s\right) =\left(\left\lfloor\left(bn\right)^{a}\right\rfloor  + t\right)^{s}\int_{0}^{\infty }x^{s - 1}e^{ - x\left(\left\lfloor\left(bn\right)^{a}\right\rfloor  + t\right)}dx
	\end{equation*}

	\begin{equation*}
		\therefore\frac{1}{\left(\left\lfloor\left(bn\right)^{a}\right\rfloor  + t\right)^{s}} =\frac{1}{\Gamma\left(s\right)}\int_{0}^{\infty }x^{s - 1}e^{ - x\left(\left\lfloor\left(bn\right)^{a}\right\rfloor  + t\right)}dx
	\end{equation*}

	\begin{equation}
		\therefore \sum_{n = 0}^{\infty }\frac{1}{\left(\left\lfloor\left(bn\right)^{a}\right\rfloor  + t\right)^{s}} =\frac{1}{t^{s}} + \sum_{n = 1}^{\infty }\frac{1}{\Gamma\left(s\right)}\int_{0}^{\infty }x^{s - 1}e^{ - x\left(\left\lfloor\left(bn\right)^{a}\right\rfloor  + t\right)}dx
	\end{equation}

	The infinite sum in the Right$-$Hand Side (RHS) of the equation \(\left(7\right)\)
	
	\begin{equation*}
		\sum_{n = 1}^{\infty }\frac{1}{\Gamma\left(s\right)}\int_{0}^{\infty }x^{s - 1}e^{ - xt}e^{ - x\left(\left\lfloor\left(bn\right)^{a}\right\rfloor\right)}dx
	\end{equation*}

	The number of repetitions of \( \left\lfloor\left(bn\right)^{a}\right\rfloor\) is given by the function \( f\left(m\right)\) (A general concept, displayed in First method, used in both methods)
	
	\begin{equation*}
		\therefore \sum_{n = 1}^{\infty }\frac{1}{\Gamma\left(s\right)}\int_{0}^{\infty }x^{s - 1}e^{ - xt}e^{ - x\left(\left\lfloor\left(bn\right)^{a}\right\rfloor\right)}dx = \sum_{n = 1}^{\infty }\frac{1}{\Gamma\left(s\right)}\int_{0}^{\infty }x^{s - 1}e^{ - xt}\left\{\left\lceil\frac{\left(n + 1\right)^{\frac{1}{a}}}{b}\right\rceil  -\left\lceil\frac{n^{\frac{1}{a}}}{b}\right\rceil\right\} e^{ - nx}dx
	\end{equation*}

	\begin{equation*}
		= \sum_{n = 1}^{\infty }\frac{\left\{\left\lceil\frac{\left(n + 1\right)^{\frac{1}{a}}}{b}\right\rceil  -\left\lceil\frac{n^{\frac{1}{a}}}{b}\right\rceil\right\} }{\Gamma\left(s\right)}\int_{0}^{\infty }x^{s - 1}e^{ -\left(n + t\right)x}dx
	\end{equation*}

	\begin{equation*}
		= \sum_{n = 1}^{\infty }\frac{\left\{\left\lceil\frac{\left(n + 1\right)^{\frac{1}{a}}}{b}\right\rceil  -\left\lceil\frac{n^{\frac{1}{a}}}{b}\right\rceil\right\} }{\Gamma\left(s\right)} \int_{0}^{\infty }\frac{\left((n + t)x\right)^{s - 1}e^{ - nx}(n + t)dx}{(n + t)^{s}}
	\end{equation*}

	\begin{equation*}
		= \sum_{n = 1}^{\infty }\frac{\left\{\left\lceil\frac{\left(n + 1\right)^{\frac{1}{a}}}{b}\right\rceil  -\left\lceil\frac{n^{\frac{1}{a}}}{b}\right\rceil\right\} }{(n + t)^{s}\ast \Gamma\left(s\right)} \int_{0}^{\infty }\left((n + t)x\right)^{s - 1}e^{ - nx}(n + t)dx
	\end{equation*}

	\begin{equation*}
		= \sum_{n = 1}^{\infty }\frac{\left\{\left\lceil\frac{\left(n + 1\right)^{\frac{1}{a}}}{b}\right\rceil  -\left\lceil\frac{n^{\frac{1}{a}}}{b}\right\rceil\right\} }{(n + t)^{s}\ast \Gamma\left(s\right)}\Gamma\left(s\right)~~~\left(~\because \int_{0}^{\infty }\left(kx\right)^{s - 1}e^{ - nx}(k)dx = \Gamma\left(s\right)~\right)
	\end{equation*}

	\begin{equation*}
		= \sum_{n = 1}^{\infty }\frac{\left\{\left\lceil\frac{\left(n + 1\right)^{\frac{1}{a}}}{b}\right\rceil  -\left\lceil\frac{n^{\frac{1}{a}}}{b}\right\rceil\right\} }{(n + t)^{s}}
	\end{equation*}

	Hence from equation (7)
	
	\begin{equation*}
		\sum_{n = 0}^{\infty }\frac{1}{\left(\left\lfloor\left(bn\right)^{a}\right\rfloor  + t\right)^{s}} =\frac{1}{t^{s}} + \sum_{n = 1}^{\infty }\frac{\left\lceil\frac{\left(n + 1\right)^{\frac{1}{a}}}{b}\right\rceil  -\left\lceil\frac{n^{\frac{1}{a}}}{b}\right\rceil }{(n + t)^{s}} = \sum_{n = 0}^{\infty }\frac{\left\lceil\frac{\left(n + 1\right)^{\frac{1}{a}}}{b}\right\rceil  -\left\lceil\frac{n^{\frac{1}{a}}}{b}\right\rceil }{(n + t)^{s}}
	\end{equation*}
	
	\textbf{SERIES II:}
	
	Again, take Gamma function
	
	\begin{equation*}
		\Gamma\left(s\right) = \int_{0}^{\infty }y^{s - 1}e^{ - y}dy
	\end{equation*}

	Take \( y = x\left(\left\lceil\left(bn\right)^{a}\right\rceil  + t\right)\)
	
	\( \Rightarrow dy = \left(\left\lceil\left(bn\right)^{a}\right\rceil  + t\right)dx\), \( x\rightarrow 0\) as \( y\rightarrow 0\) and \( x\rightarrow \infty\) as \( y\rightarrow \infty\)
	
	\begin{equation*}
		\therefore  \Gamma\left(s\right) = \int_{0}^{\infty }\left(x\left(\left\lceil\left(bn\right)^{a}\right\rceil  + t\right)\right)^{s - 1}e^{ - x\left(\left\lceil\left(bn\right)^{a}\right\rceil  + t\right)}\left(\left\lceil\left(bn\right)^{a}\right\rceil  + t\right)dx
	\end{equation*}

	\begin{equation*}
		\therefore  \Gamma\left(s\right) =\left(\left\lceil\left(bn\right)^{a}\right\rceil  + t\right)^{s}\int_{0}^{\infty }x^{s - 1}e^{ - x\left(\left\lceil\left(bn\right)^{a}\right\rceil  + t\right)}dx
	\end{equation*}

	\begin{equation*}
		\therefore\frac{1}{\left(\left\lceil\left(bn\right)^{a}\right\rceil  + t\right)^{s}} =\frac{1}{\Gamma\left(s\right)}\int_{0}^{\infty }x^{s - 1}e^{ - x\left(\left\lceil\left(bn\right)^{a}\right\rceil  + t\right)}dx
	\end{equation*}

	\begin{equation}
		\therefore \sum_{n = 0}^{\infty }\frac{1}{\left(\left\lceil\left(bn\right)^{a}\right\rceil  + t\right)^{s}} =\frac{1}{t^{s}} + \sum_{n = 1}^{\infty }\frac{1}{\Gamma\left(s\right)}\int_{0}^{\infty }x^{s - 1}e^{ - x\left(\left\lceil\left(bn\right)^{a}\right\rceil  + t\right)}dx                                                                      
	\end{equation}

	The infinite sum in the Right$-$Hand Side (RHS) of the equation \(\left(8\right)\)
	
	\begin{equation*}
		\sum_{n = 1}^{\infty }\frac{1}{\Gamma\left(s\right)}\int_{0}^{\infty }x^{s - 1}e^{ - xt}e^{ - x\left\lceil\left(bn\right)^{a}\right\rceil }dx
	\end{equation*}

	The number of repetitions of \(\left\lceil\left(bn\right)^{a}\right\rceil\) is given by the function \( g\left(k\right)\) (A general concept, displayed in First method, used in both methods)
	
	\begin{equation*}
		\therefore \sum_{n = 1}^{\infty }\frac{1}{\Gamma\left(s\right)}\int_{0}^{\infty }x^{s - 1}e^{ - xt}e^{ - x\left(\left\lceil\left(bn\right)^{a}\right\rceil\right)}dx = \sum_{n = 1}^{\infty }\frac{1}{\Gamma\left(s\right)}\int_{0}^{\infty }x^{s - 1}e^{ - xt}\left\{\left\lfloor\frac{n^{\frac{1}{a}}}{b}\right\rfloor  -\left\lfloor\frac{(n - 1)^{\frac{1}{a}}}{b}\right\rfloor\right\} e^{ - nx}dx
	\end{equation*}

	\begin{equation*}
		= \sum_{n = 1}^{\infty }\frac{\left\{\left\lfloor\frac{n^{\frac{1}{a}}}{b}\right\rfloor  -\left\lfloor\frac{(n - 1)^{\frac{1}{a}}}{b}\right\rfloor\right\} }{\Gamma\left(s\right)}\int_{0}^{\infty }x^{s - 1}e^{ -\left(n + t\right)x}dx
	\end{equation*}

	\begin{equation*}
		= \sum_{n = 1}^{\infty }\frac{\left\{\left\lfloor\frac{n^{\frac{1}{a}}}{b}\right\rfloor  -\left\lfloor\frac{(n - 1)^{\frac{1}{a}}}{b}\right\rfloor\right\} }{\Gamma\left(s\right)} \int_{0}^{\infty }\frac{\left((n + t)x\right)^{s - 1}e^{ - nx}(n + t)dx}{(n + t)^{s}}
	\end{equation*}

	\begin{equation*}
		= \sum_{n = 1}^{\infty }\frac{\left\{\left\lfloor\frac{n^{\frac{1}{a}}}{b}\right\rfloor  -\left\lfloor\frac{(n - 1)^{\frac{1}{a}}}{b}\right\rfloor\right\} }{(n + t)^{s}\ast \Gamma\left(s\right)} \int_{0}^{\infty }\left((n + t)x\right)^{s - 1}e^{ - nx}(n + t)dx
	\end{equation*}

	\begin{equation*}
		= \sum_{n = 1}^{\infty }\frac{\left\{\left\lfloor\frac{n^{\frac{1}{a}}}{b}\right\rfloor  -\left\lfloor\frac{(n - 1)^{\frac{1}{a}}}{b}\right\rfloor\right\} }{(n + t)^{s}\ast \Gamma\left(s\right)}\Gamma\left(s\right)                  \left(\because \int_{0}^{\infty }\left(kx\right)^{s - 1}e^{ - nx}(k)dx = \Gamma\left(s\right)\right)
	\end{equation*}

	\begin{equation*}
		= \sum_{n = 1}^{\infty }\frac{\left\{\left\lfloor\frac{n^{\frac{1}{a}}}{b}\right\rfloor  -\left\lfloor\frac{(n - 1)^{\frac{1}{a}}}{b}\right\rfloor\right\} }{(n + t)^{s}}
	\end{equation*}

	Hence from equation (8)
	\begin{equation*}
		\sum_{n = 0}^{\infty }\frac{1}{\left(\left\lceil\left(bn\right)^{a}\right\rceil  + t\right)^{s}} =\frac{1}{t^{s}} + \sum_{n = 1}^{\infty }\frac{\left\lfloor\frac{n^{\frac{1}{a}}}{b}\right\rfloor  -\left\lfloor\frac{(n - 1)^{\frac{1}{a}}}{b}\right\rfloor }{(n + t)^{s}}
	\end{equation*}
	The Right hand sides in the equations \((7)\) and  \((8)\) are observed to be faster convergent. (See open problem 4.2) 
	
	\subsection{Deductions:}
	
	\textbf{DEDUCTION 4:} Taking \( a =\frac{1}{q} ; q\in \mathbb{N}, Re\left(s\right)>q\) and \( b  =  1\):
	
	\begin{equation*}
		_{  }^{F}\zeta_{1}^{\frac{1}{q}}\left(s,t\right) = \sum_{n = 0}^{\infty }\frac{1}{\left(\left\lfloor n^{\frac{1}{q}}\right\rfloor  + t\right)^{s}} = \sum_{m = 0}^{q - 1}\sum_{k = m}^{q - 1}\begin{pmatrix}
			q \\ 
			m \\ 
		\end{pmatrix}\begin{pmatrix}
			q \\ 
			q - k - 1 \\ 
		\end{pmatrix}\left( - t\right)^{q - k - 1}\zeta\left(s - m,t\right)
	\end{equation*}
	
	\begin{equation*}
		_{  }^{F}\zeta_{1}^{\frac{1}{q}}\left(s,t\right) = \int_{0}^{\infty }\frac{P\left(x,s,t,q\right)}{1-e^{-x}}dx. ~~~~~ where ~~ P\left(x,s,t,q\right) =\sum_{m = 0}^{q - 1}\sum_{k = m}^{q - 1}\begin{pmatrix}
			q \\ 
			m \\ 
		\end{pmatrix}\begin{pmatrix}
			q \\ 
			q - k - 1 \\ 
		\end{pmatrix}\left( - t\right)^{q - k - 1}\left(\Gamma(s-m)\right)^{-1}e^{-tx}
	\end{equation*}
	
	\textbf{DEDUCTION 5:} Taking \( a =\frac{1}{q}; q\in \mathbb{N}, Re\left(s\right)>q\) and \( b  =  1\):
	
	\begin{equation*}
		_{  }^{C}\zeta_{1}^{\frac{1}{q}}\left(s,t\right) = \sum_{n = 0}^{\infty }\frac{1}{\left(\left \lceil n^{\frac{1}{q}}\right\rceil  + t\right)^{s}} = \sum_{m = 0}^{q - 1}\sum_{k = m}^{q - 1}\left( - 1\right)^{q - m + 1}\begin{pmatrix}
			q \\ 
			m \\ 
		\end{pmatrix}\begin{pmatrix}
			q \\ 
			q - k - 1 \\ 
		\end{pmatrix}\left(t\right)^{q - k - 1} \zeta\left(s - m,t\right)
	\end{equation*}
	
	\begin{equation*}
		_{  }^{F}\zeta_{1}^{\frac{1}{q}}\left(s,t\right) = \int_{0}^{\infty }\frac{Q\left(x,s,t,q\right)}{1-e^{-x}}dx. ~~~~~ where ~~ Q\left(x,s,t,q\right) \sum_{m = 0}^{q - 1}\sum_{k = m}^{q - 1}\left( - 1\right)^{q - m + 1}\begin{pmatrix}
			q \\ 
			m \\ 
		\end{pmatrix}\begin{pmatrix}
			q \\ 
			q - k - 1 \\ 
		\end{pmatrix}\left(t\right)^{q - k - 1} \left(\Gamma(s-m)\right)^{-1}e^{-tx}
\end{equation*}
	
	\textbf{DEDUCTION 6: }Taking \( a =\frac{1}{q}; q\in \mathbb{N},Re\left(s\right) > q, b  =  1\) and \( t = 1\):
	
	\begin{equation*}
		_{  }^{F}\zeta_{1}^{\frac{1}{q}}\left(s,1\right) =_{  }^{F}\zeta_{1}^{\frac{1}{q}}\left(s\right)  = \sum_{n = 1}^{\infty }\frac{1}{\left\lfloor n^{\frac{1}{q}}\right\rfloor ^{s}} = \sum_{m = 0}^{q - 1}\begin{pmatrix}
			q \\ 
			m \\ 
		\end{pmatrix}\zeta\left(s - m\right)
	\end{equation*}

	Using the definition of the Zeta function, this can also be re$-$written in the integral form as:
	
	\begin{equation*}
		_{  }^{F}\zeta_{1}^{\frac{1}{q}}\left(s\right) = \int_{0}^{\infty }\frac{P\left(x,s,q\right)}{e^{x} - 1}dx. ~~~~~ where ~~ P\left(x,s,q\right) = \sum_{t = 0}^{q - 1}\left(\Gamma\left(s - t\right)\right)^{ - 1}\binom{q}{t}x^{s - t - 1}
	\end{equation*}
	
	Furthermore, given the analytic continuation of the Zeta function we get that $ _{  }^{F}\zeta_{1}^{\frac{1}{q}}\left(s\right) $ can be defined even for $ Re\left(s\right) < q $ and hence $ _{  }^{F}\zeta_{1}^{\frac{1}{q}}\left(s\right) $ can be defined for $ \forall s \in \mathbb{C} $ such that $ Re(s) \neq q $
		
	\textbf{DEDUCTION 7: }Taking \( ,a =\frac{1}{q}; q\in \mathbb{N},Re\left(s\right) > q, b  =  1\) and \( t = 1\):
	\begin{align*}
		&_{  }^{C}\zeta_{1}^{\frac{1}{q}}\left(s,1\right) = _{  }^{C}\zeta_{1}^{\frac{1}{q}}\left(s\right)= \sum_{n = 1}^{\infty }\frac{1}{\left\lceil n^{\frac{1}{q}}\right\rceil ^{s}} = \sum_{m = 0}^{q - 1}\left( - 1\right)^{q - m + 1}\begin{pmatrix}
			q \\ 
			m \\ 
		\end{pmatrix}\zeta\left(s - m\right) \\
		 \therefore ~ & _{  }^{C}\zeta_{1}^{\frac{1}{q}}\left(s\right) = \int_{0}^{\infty }\frac{Q\left(x,s,q\right)}{e^{x} - 1}dx. ~~~~~ where ~ Q\left(x,s,q\right) = \sum_{t = 0}^{q - 1}\left( - 1\right)^{q - t + 1}\left(\Gamma\left(s - t\right)\right)^{ - 1}\binom{q}{t}x^{s - t - 1}
	\end{align*}
	Again, given the analytic continuation of the Zeta function we get that $ _{  }^{C}\zeta_{1}^{\frac{1}{q}}\left(s\right) $ can be defined even for $ Re\left(s\right) < q $ and hence $ _{  }^{C}\zeta_{1}^{\frac{1}{q}}\left(s\right) $ can be defined for $ \forall s \in \mathbb{C} $ such that $ Re(s) \neq q $

	\vspace{1\baselineskip}
	All series in deductions 4,5, 6 and 7 have poles at \( Re\left(s\right) = q\), but following can be simply observed for \( q\neq 1\):
	
	\begin{equation}
		_{  }^{C}\zeta_{1}^{\frac{1}{q}}\left(q,t\right) -_{  }^{F}\zeta_{1}^{\frac{1}{q}}\left(q,t\right) = \sum_{m = 0}^{q - 2}\sum_{k = m}^{q - 2}\begin{pmatrix}
			q \\ 
			m \\ 
		\end{pmatrix}\begin{pmatrix}
			q \\ 
			q - k - 1 \\ 
		\end{pmatrix}\zeta\left(q - m,t\right)\left\{\left( - t\right)^{q - k - 1} +\left( - 1\right)^{q - m}\left(t\right)^{q - k - 1}\right\} 
	\end{equation}
	
	\begin{equation}
		_{  }^{C}\zeta_{1}^{\frac{1}{q}}\left(q\right) -_{  }^{F}\zeta_{1}^{\frac{1}{q}}\left(q\right)~ = ~ \sum_{n = 1}^{\infty }\left(\frac{1}{\left\lfloor n^{\frac{1}{q}}\right\rfloor^{q}} -\frac{1}{\left\lceil n^{\frac{1}{q}}\right\rceil^{q}}\right) ~ = ~ \sum_{t = 0}^{q - 2}\begin{pmatrix}
			q \\ 
			t \\ 
		\end{pmatrix} \zeta\left(q - t\right)\left[1 +\left( - 1\right)^{q - t}\right]
	\end{equation}

	This shows that even if the set two series may individually have poles at \( s = q\) but their difference is convergent.
	
	\subsection{Results for Specific Values}
	\textbf{(I):} Taking \(  a =  b = 1\) in series 1 and 2, both reduce to Hurwitz$-$Zeta function. (Consider series 1)
	
	\begin{equation*}
		_{  }^{F}\zeta_{1}^{1}\left(s,t\right) = \sum_{n = 0}^{\infty }\frac{1}{\left\lfloor\left(1\ast n\right)^{1} + t\right\rfloor^{s}} = \sum_{n = 0}^{\infty }\frac{1}{\left(n + t\right)^{s}} = \zeta\left(s,t\right)
	\end{equation*}
	
	\textbf{(II):} Taking \(  a =  b = t = 1\) in series 1 and 2, both reduce to Riemann$-$Zeta function. (Consider series 2)
	
	\begin{equation*}
		_{  }^{C}\zeta_{1}^{1}\left(s\right) = \sum_{n = 0}^{\infty }\frac{1}{\left(\left\lceil\left(1\ast n\right)^{1}\right\rceil  + 1\right)^{s}} = \sum_{n = 1}^{\infty }\frac{1}{n^{s}} = \zeta\left(s\right)
	\end{equation*}

	Following are special cases of series 1 and 2 which are observed solely by intuition:
	
	\begingroup
	\large
	\doublespacing
	\begin{tabular}{l l}
	$\left(1\right): ~ _{}^{F}\zeta_{2}^{\frac{1}{2}}\left(s\right) = \sum\limits_{n = 1}^{\infty }\frac{2\left\lfloor\frac{n}{2}\right\rfloor  + 1}{n^{s}}$ & $\left(2\right): ~_{  }^{C}\zeta_{2}^{\frac{1}{2}}\left(s\right) = \sum\limits_{n = 1}^{\infty }\frac{2\left\lfloor\frac{n}{2}\right\rfloor }{n^{s}}$\\
	$ \left(3\right):~_{  }^{F}\zeta_{3}^{\frac{1}{2}}\left(s\right) = \sum\limits_{n = 1}^{\infty }\frac{2\left\lfloor\frac{n}{3}\right\rfloor  + 1}{n^{s}} $ & $ \left(4\right):~_{  }^{C}\zeta_{3}^{\frac{1}{2}}\left(s\right) = \sum\limits_{n = 1}^{\infty }\frac{\left\lfloor\frac{2n}{3}\right\rfloor }{n^{s}} $ \\ $
		\left(5\right):~_{  }^{F}\zeta_{4}^{\frac{1}{2}}\left(s\right) = \sum\limits_{n = 1}^{\infty }\frac{1}{\left\lfloor 2\sqrt{n}\right\rfloor^{s}} = \sum\limits_{n = 1}^{\infty }\frac{\left\lceil\frac{n}{2}\right\rceil  +\left( - 1\right)^{n}}{n^{s}} $  & $ \left(6\right):~_{  }^{C}\zeta_{4}^{\frac{1}{2}}\left(s\right) = \sum\limits_{n = 1}^{\infty }\frac{1}{\left\lceil 2\sqrt{n}\right\rceil^{s}} = \sum\limits_{n = 1}^{\infty }\frac{\left\lfloor\frac{n}{2}\right\rfloor }{n^{s}} $ \\
	 	$ \left(7\right):~_{  }^{F}\zeta_{2}^{\frac{1}{3}}\left(s\right) = \sum\limits_{n = 1}^{\infty }\frac{\frac{3n\left(n + 1\right)}{2} +\left\lfloor\frac{n}{2}\right\rfloor  -\left\lfloor\frac{n - 1}{2}\right\rfloor }{n^{s}} $ & $ \left(8\right):~_{  }^{C}\zeta_{2}^{\frac{1}{3}}\left(s\right) = \sum\limits_{n = 1}^{\infty }\frac{\frac{3n\left(n - 1\right)}{2} +\left\lfloor\frac{n}{2}\right\rfloor  -\left\lfloor\frac{n - 1}{2}\right\rfloor }{n^{s}}$ \\ 
		$ \left(9\right):~_{  }^{F}\zeta_{3}^{\frac{1}{3}}\left(s\right) = \sum\limits_{n = 1}^{\infty }\frac{n\left(n + 1\right) +\left\lfloor\frac{n}{3}\right\rfloor  -\left\lfloor\frac{n - 1}{3}\right\rfloor }{n^{s}} $ & $ \left(10\right):~_{  }^{C}\zeta_{3}^{\frac{1}{3}}\left(s\right) = \sum\limits_{n = 1}^{\infty }\frac{n\left(n - 1\right) +\left\lfloor\frac{n}{3}\right\rfloor  -\left\lfloor\frac{n - 1}{3}\right\rfloor }{n^{s}} $\\
		$ \left(11\right):~_{  }^{F}\zeta_{4}^{\frac{1}{3}}\left(s\right) = \sum\limits_{n = 1}^{\infty }\frac{\left\lfloor\frac{3n\left(n + 1\right)}{4}\right\rfloor  +\left\lfloor\frac{n}{2}\right\rfloor  -\left\lfloor\frac{n - 1}{2}\right\rfloor }{n^{s}}$ & $ \left(12\right):~_{  }^{F}\zeta_{4}^{\frac{1}{3}}\left(s\right) = \sum\limits_{n = 1}^{\infty }\frac{\left\lfloor\frac{3n\left(n - 1\right)}{4}\right\rfloor  +\left\lfloor\frac{n}{2}\right\rfloor  -\left\lfloor\frac{n - 1}{2}\right\rfloor }{n^{s}} $ \\
	 	$ \left(13\right):~_{  }^{F}\zeta_{2}^{\frac{1}{4}}\left(s\right) = \sum\limits_{n = 1}^{\infty }\frac{n\left(n + 1\right)\left(2n + 1\right) + 2\left\lfloor\frac{n}{2}\right\rfloor  + 1}{n^{s}} $ & 
	 	$  \left(14\right):~_{  }^{F}\zeta_{2}^{\frac{1}{4}}\left(s\right) = \sum\limits_{n = 1}^{\infty }\frac{\left(n - 1\right)n\left(2n - 1\right) + 2\left\lfloor\frac{n}{2}\right\rfloor }{n^{s}} $ \\
	 	$ \left(15\right):~_{  }^{F}\zeta_{2}^{\frac{1}{5}}\left(s\right) = \sum\limits_{n = 1}^{\infty }\frac{\frac{5}{2}n\left(n + 1\right)\left(n^{2} + n + 1\right) +\left\lfloor\frac{n}{2}\right\rfloor  -\left\lfloor\frac{n - 1}{2}\right\rfloor }{n^{s}}$ & \\
	 	$ \left(16\right):~_{  }^{C}\zeta_{2}^{\frac{1}{5}}\left(s\right) = \sum\limits_{n = 1}^{\infty }\frac{\frac{5}{2}\left(n - 1\right)n\left(n^{2} - n + 1\right) +\left\lfloor\frac{n}{2}\right\rfloor  -\left\lfloor\frac{n - 1}{2}\right\rfloor }{n^{s}} $ & \\
	\end{tabular}
	\endgroup \vspace{2mm}
	
	One can go further with \( a =\frac{1}{6},\frac{1}{7},\frac{1}{8},\ldots\) with different values of \( b\) for both functions.
	
	\section{Conclusion and Open Problems} \label {section 5}
		
		This paper introduces a set of a new type of partial summation formulas and generalisation of Hurwitz Zeta function using floor and ceiling functions. Very fundamental yet powerful principle of mathematical induction is used to prove the formulas. Two functions \( f\left(m\right)\) and \( g(k)\) are used to get the number of repetitions of \(\left\lfloor\left(bn\right)^{a}\right\rfloor\) and \(\left\lceil\left(bn\right)^{a}\right\rceil\) for consecutive natural numbers \( n\) and are used in getting the equivalent series for \( \sum\limits_{n = 1}^{\infty }\left\lfloor\left(bn\right)^{a}\right\rfloor^{ - s}\) and \( \sum\limits_{n = 1}^{\infty }\left\lceil\left(bn\right)^{a}\right\rceil^{ - s}\). It is shown that both F$-$Hurwitz and C$-$Hurwitz Zeta functions reduce to Hurwitz$-$Zeta function and Riemann$-$Zeta function for particular values. A set of results and special cases helps to understand the behaviour of the formulas and series at specific values. The extensions done in the paper may be helpful in improving the measurement quality as well expanding the domain of research.
		
		\vspace{1\baselineskip}
		\textbf{Problem 4.1: }Finding proofs for the explicit formulas for \(\left(i\right)\sum\limits_{i = 1}^{n}\left\lfloor i^{a}\right\rfloor^{p} \& \left(ii\right)\sum\limits_{i = 1}^{n}\left\lceil i^{a}\right\rceil^{p} , p\in \mathbb{N}, a\in \mathbb{R}^{ + }\) (Section \ref{section 2}).
		
		\vspace{2mm}
		 		
	\textbf{Problem 4.2(A): }Are the following assumptions true?
	
	If \( Re\left(s\right)>q\) then
	
	\textbf{(I)}
	
	\begin{align*}
		_{  }^{F}\zeta_{2}^{\frac{1}{q}}\left(s\right) = & \sum_{n = 1}^{\infty }\frac{1}{\left\lfloor\left(2n\right)^{\frac{1}{q}}\right\rfloor^{s}} = \sum\limits_{m = 1}^{\infty }\frac{\frac{q\left(q - 1\right)}{2}\left\{ \sum\limits_{t =\left\lceil\frac{q}{2}\right\rceil  -\left\lfloor\frac{q}{2}\right\rfloor }^{q - 2}\left(\sum\limits_{i = 1}^{m}\left(\frac{i^{t}\left(1 +\left( - 1\right)^{q - t}\right)}{2}\right)\right)\right\}  + y_{m}}{m^{s}} \\
		& y_{m} =\left\{ 
		\begin{array}{c}
			2\left\lfloor\frac{m}{2}\right\rfloor  + 1 , q even \\ 
			\left\lfloor\frac{m}{2}\right\rfloor  -\left\lfloor\frac{m - 1}{2}\right\rfloor , q odd \\ 
		\end{array}\right.
	\end{align*}

	\vspace{1\baselineskip}
	\textbf{(II)}
	
	\begin{align*}
		_{  }^{C}\zeta_{2}^{\frac{1}{q}}\left(s\right) = & \sum_{n = 1}^{\infty }\frac{1}{\left\lceil\left(2n\right)^{\frac{1}{q}}\right\rceil^{s}} = \sum\limits_{m = 1}^{\infty }\frac{\frac{q\left(q - 1\right)}{2}\left\{ \sum\limits_{t =\left\lceil\frac{q}{2}\right\rceil  -\left\lfloor\frac{q}{2}\right\rfloor }^{q - 2}\left(\sum\limits_{i = 1}^{m - 1}\left(\frac{i^{t}\left(1 +\left( - 1\right)^{q - t}\right)}{2}\right)\right)\right\}  + z_{m}}{m^{s}} \\
		& z_{m} =\left\{ 
		\begin{array}{c}
			2\left\lfloor\frac{m}{2}\right\rfloor  , q even \\ 
			\left\lfloor\frac{m}{2}\right\rfloor  -\left\lfloor\frac{m - 1}{2}\right\rfloor , q odd \\ 
		\end{array}\right.
	\end{align*}

	\textbf{(B): }If these assumption holds true, can one go on to find a general formula for the equivalents of \(_{  }^{F}\zeta_{k}^{\frac{1}{q}}\left(s\right) \& _{  }^{C}\zeta_{k}^{\frac{1}{q}}\left(s\right) ,Re\left(s\right)>q , k\in \mathbb{N}\textbackslash\left\{ 1,2\right\}\)? \vspace{2mm}
	
	\textbf{Problem 4.3: }In the subsection \ref{subsection 4.1}, for both C$-$Hurwitz and F$-$Hurwitz Zeta functions, the corresponding infinite series on the right$-$hand side (equivalent series) converge faster to the exact values of the Zeta functions, one can prove the same mathematically. \vspace{2mm}
	
	\textbf{Problem 4.4: } Consider the following double series
	
	\begin{equation*}
		S_{n} = \sum_{i = 1}^{n}\sum_{j = 1}^{\frac{i^{p}}{j^{q`}}}a_{i,j}
	\end{equation*}
	
	If \(  a_{i,j} = 1\) then does it hold true that the upper bound of time complexity (Big$-$O) for the given series is \(  O\left(\left\lfloor n^{\frac{p + q + 1}{q + 1}}\right\rfloor\right)\). \vspace{2mm}
	
	\textbf{Problem 4.5: }The Riemann hypothesis is a very well$-$known unsolved problem, involves the Zeta function. Both new Zeta functions reduce to the Riemann$-$Zeta function for specific values. Can these new Zeta functions or their deduction be useful in further understanding the hypothesis? \vspace{2mm}
	
	\vspace{1\baselineskip}
	\begin{center}
		{\large A}CKNOWLEDGEMENT
	\end{center}

	The authors would like to acknowledge Dr. Nishant Doshi of Pandit Deendayal Energy University, Gandhinagar, Gujarat, India, for bringing us a very interesting exercise, from which we were able to work on all the results. Authors would like to thank Kenneth Beitler for contacting us and for giving a very helpful suggestion for proving the results. We are especially grateful to have guidance of Dr. Anand Sengupta and Dr. Atul Dixit of Indian Institute of Technology, Gandhinagar by the means of video conference and discussions. Also, we are thankful to Ms. Meghna Parikh and Mr. Bhashin Thakore for valuable discussion throughout the preparation.


\begin{thebibliography}{10}
		
		\bibitem{Apostol76}
		Tom Apostol,
		\textit{Introduction to Analytic Number Theory},
		Springer, 
		1976.
		
		\bibitem{CG96}
		J. H. Conway and R. Guy, ,
		\textit{The Book of Numbers},
		Springer, 
		1996.
		
		\bibitem{MNR17}
		H. Montgomery, A. Nikeghbali, M. Rassias, 
		\textit{Exploring the Riemann Zeta Function: 190 years from Riemann’s Birth},
		Springer, 
		2017.
		
		\bibitem{HW80}
		G. H. Hardy and E. M. Wright,
		\textit{An introduction to the theory of numbers},
		Fifth ed.,
		Clarendon Press, Oxford, 
		1980.
		
		\bibitem{GKP89}
		R. L. Graham, D. E. Knuth and O. Patashnik,
		\textit{Concrete mathematics: A Foundation for Computer Science},
		Addison-Wesley,
		1994.
	\end{thebibliography}
\end{document}